\documentclass[12pt]{article}
\usepackage{amsmath,amssymb,amsbsy,amsfonts,amsthm,latexsym,
                     amsopn,amstext,amsxtra,euscript,amscd}

\def\bbbc{{\mathchoice {\setbox0=\hbox{$\displaystyle\rm C$}\hbox{\hbox
to0pt{\kern0.4\wd0\vrule height0.9\ht0\hss}\box0}}
{\setbox0=\hbox{$\textstyle\rm C$}\hbox{\hbox
to0pt{\kern0.4\wd0\vrule height0.9\ht0\hss}\box0}}
{\setbox0=\hbox{$\scriptstyle\rm C$}\hbox{\hbox
to0pt{\kern0.4\wd0\vrule height0.9\ht0\hss}\box0}}
{\setbox0=\hbox{$\scriptscriptstyle\rm C$}\hbox{\hbox
to0pt{\kern0.4\wd0\vrule height0.9\ht0\hss}\box0}}}}
\def\bbbq{{\mathchoice {\setbox0=\hbox{$\displaystyle\rm
Q$}\hbox{\raise
0.15\ht0\hbox to0pt{\kern0.4\wd0\vrule height0.8\ht0\hss}\box0}}
{\setbox0=\hbox{$\textstyle\rm Q$}\hbox{\raise
0.15\ht0\hbox to0pt{\kern0.4\wd0\vrule height0.8\ht0\hss}\box0}}
{\setbox0=\hbox{$\scriptstyle\rm Q$}\hbox{\raise
0.15\ht0\hbox to0pt{\kern0.4\wd0\vrule height0.7\ht0\hss}\box0}}
{\setbox0=\hbox{$\scriptscriptstyle\rm Q$}\hbox{\raise
0.15\ht0\hbox to0pt{\kern0.4\wd0\vrule height0.7\ht0\hss}\box0}}}}
\def\bbbt{{\mathchoice {\setbox0=\hbox{$\displaystyle\rm
T$}\hbox{\hbox to0pt{\kern0.3\wd0\vrule height0.9\ht0\hss}\box0}}
{\setbox0=\hbox{$\textstyle\rm T$}\hbox{\hbox
to0pt{\kern0.3\wd0\vrule height0.9\ht0\hss}\box0}}
{\setbox0=\hbox{$\scriptstyle\rm T$}\hbox{\hbox
to0pt{\kern0.3\wd0\vrule height0.9\ht0\hss}\box0}}
{\setbox0=\hbox{$\scriptscriptstyle\rm T$}\hbox{\hbox
to0pt{\kern0.3\wd0\vrule height0.9\ht0\hss}\box0}}}}
\def\bbbs{{\mathchoice
{\setbox0=\hbox{$\displaystyle     \rm S$}\hbox{\raise0.5\ht0\hbox
to0pt{\kern0.35\wd0\vrule height0.45\ht0\hss}\hbox
to0pt{\kern0.55\wd0\vrule height0.5\ht0\hss}\box0}}
{\setbox0=\hbox{$\textstyle        \rm S$}\hbox{\raise0.5\ht0\hbox
to0pt{\kern0.35\wd0\vrule height0.45\ht0\hss}\hbox
to0pt{\kern0.55\wd0\vrule height0.5\ht0\hss}\box0}}
{\setbox0=\hbox{$\scriptstyle      \rm S$}\hbox{\raise0.5\ht0\hbox
to0pt{\kern0.35\wd0\vrule height0.45\ht0\hss}\raise0.05\ht0\hbox
to0pt{\kern0.5\wd0\vrule height0.45\ht0\hss}\box0}}
{\setbox0=\hbox{$\scriptscriptstyle\rm S$}\hbox{\raise0.5\ht0\hbox
to0pt{\kern0.4\wd0\vrule height0.45\ht0\hss}\raise0.05\ht0\hbox
to0pt{\kern0.55\wd0\vrule height0.45\ht0\hss}\box0}}}}
\def\bbbz{{\mathchoice {\hbox{$\sf\textstyle Z\kern-0.4em Z$}}
{\hbox{$\sf\textstyle Z\kern-0.4em Z$}}
{\hbox{$\sf\scriptstyle Z\kern-0.3em Z$}}
{\hbox{$\sf\scriptscriptstyle Z\kern-0.2em Z$}}}}

\def\cA{\mathcal A}
\def\cF{\mathcal F}

\def\cB{\mathcal B}

\def\cN{\mathcal N}
\def\cM{\mathcal M}
\def\cS{\mathcal S}
\def\cU{\mathcal U}

\def\cX{\mathcal X}

\def\Z{\mathbb{Z}}

\def\R{\mathbb{R}}
\def\T{\mathbb{T}}

\renewcommand{\vec}[1]{\mathbf{#1}}

\def\vec#1{\mathbf{#1}}
\def\inv#1{\overline{#1}}

\def\bbbc{{\mathrm I\!C}}

\def\bbbq{{\mathbb Q}}

\def \Z{{\bbbz}}

\def\dist{\mathrm{dist}}

\def\Nab{\cN(\vec{a},\vec{b},q)}
\def\Hab{H(\vec{a},\vec{b},q)}
\def\NAB{\cN_\varOmega(\vec{a},\vec{b},q)}
 
\def\MAB{\cM_\varTheta(\vec{a},\vec{b},q)}
 \def\MABB{\cM_\varSigma(\vec{a},\vec{b},q)}
 
\def\NABp{\cN_\varOmega(\vec{a},\vec{b},p)}
 
\def\MABp{\cM_\varTheta(\vec{a},\vec{b},p)}
 \def\MABBp{\cM_\varSigma(\vec{a},\vec{b},p)}

\newtheorem{thm}{Theorem}
\newtheorem{lem}[thm]{Lemma}

\begin{document}

\def\\{\cr}
\def\({\left(}
\def\){\right)}
\def\fl#1{\left\lfloor#1\right\rfloor}
\def\rf#1{\left\lceil#1\right\rceil}

\title{On a Generalisation of a Lehmer Problem}

\author{ 
{\sc Igor E.~Shparlinski} \\
{Department of Computing}\\
{Macquarie University} \\
{Sydney, NSW 2109, Australia} \\
{\tt igor@ics.mq.edu.au}}

%

        \maketitle

\begin{abstract}
We consider a generalisation of the classical Lehmer problem
about the distribution of modular inverses in arithmetic 
progression, introduced by E.~Alkan, F.~Stan and A.~Zaharescu. 
Using bounds of sums of multiplicative characters instead of 
traditionally applied to this kind of problem Kloosterman sums,
we improve their results in several directions.
\end{abstract}


\section{Introduction} 

Given modulus $q\ge 2$, we denote by $\cU_q$ 
the set 
$$
\cU_q = \{n \ : \ 1 \le n < q, \ \gcd(n,q)=1\}.
$$
that is, $\# \cU_q = \varphi(q)$, the Euler function.

For $n \in \cU_q$ we use $\inv{n}$ to denote 
the modular inverse
of 
$n$, that is, $n \inv{n} \equiv 1 \pmod q$, $\inv{n}\in \cU_q$.

The classical question of D.~H.~Lehmer 
(see~\cite[Problem~F12]{Guy}) about the 
joint distribution of the parity of $n$ and $\inv{n}$
has been solved by W.~Zhang~\cite{Zha1,Zha2}.

Recently this question  has been generalised by 
E.~Alkan, F.~Stan and A.~Zaharescu~\cite{ASZ} as follows. 
Given vector 
$\vec{a} = (a_1, \ldots, a_{k+1})\in \cU_q^{k+1}$ 
 and  
$\vec{b} = (b_1, \ldots, b_{k+1})\in \Z^{k+1}$ we consider the set 
$\Nab$ of vectors  $\vec{n} = (n_1, \ldots, n_k)\in \cU_q^k$
such that 
\begin{eqnarray*}
\label{eq:Cong}
n_i &\equiv& b_i \pmod {a_i}, \qquad i =1, \ldots, k,\\
\inv{n_1 \ldots n_k} &\equiv& b_{k+1} \pmod {a_{k+1}}.
\end{eqnarray*}

Generalising several previous results of various authors, 
(for instance, of~\cite{CoZa,Zha1,Zha2,Zha3}), 
E.~Alkan, F.~Stan and A.~Zaharescu~\cite{ASZ} have 
shown that for any fixed $k$, the bound 
\begin{equation}
\label{eq:Old Nab}
\# \Nab = \frac{\varphi(q)^k}{a_1 \ldots a_{k+1}}
+ O(q^{k-1/2 + o(1)})
\end{equation} 
holds uniformly over all 
vectors $\vec{a} \in \cU_q^{k+1}$ and  $\vec{b} \in
\Z^{k+1}$.
In particular, since
\begin{equation}
\label{eq:phi}
\varphi(q) \ge c\frac{q}{\log \log (q+2)}
\end{equation}
for an absolute constant $c >0$, see~\cite[Section~I.5.4]{Ten},
we see that the bound~\eqref{eq:Old Nab} is nontrivial 
for 
\begin{equation}
\label{eq:Old Threshold}
a_1 \ldots a_{k+1} \le q^{1/2 - \delta}
\end{equation} 
for any fixed $\delta>0$,
provided that $q$ is large enough. 

The main tool of~\cite{ASZ} are bounds of exponential
sum, in particular, multidimensional Kloosterman
sums. 
Here we show that  
using bounds of multiplicative character sums,
such   the classical Polya-Vinogradov and Burgess bounds, 
see~\cite[Theorems~12.5 and 12.6]{IwKow}, 
one can improve the bound~\eqref{eq:Old Nab}.
For example,  we obtain a bound, 
which in particular implies that 
\begin{equation}
\label{eq:New Nab}
\# \Nab = \frac{\varphi(q)^k}{a_1 \ldots a_{k+1}}
+  O\(\frac{\|\vec{a}\| q^{k-1 + o(1)}}{a_1 \ldots a_{k+1}} +\frac{q^{(k+1)/2 +
o(1)}}{(a_1\ldots a_{k+1})^{1/(k+1)}} \), 
\end{equation} 
where $ \|\vec{a}\| $ is the Euclidean norm of $\vec{a}$,  
which is equivalent to~\eqref{eq:Old Nab} for  $k=2$ and 
$\|\vec{a}\| = O(1)$ and always improves it if either 
$\|\vec{a}\|$ grows together with $q$ or if
$k \ge 3$
(in this case, with respect to both dependence on $q$ and $\vec{a}$).
 
We note  that instead of~\eqref{eq:Old Threshold}, 
the bound of our Theorem~\ref{thm:main-1}  is 
nontrivial when simultaneously 
\begin{equation}
\label{eq:New Threshold}
\|\vec{a}\| \le q^{1 - \delta}
\quad \text{and}
\quad     a_1 \ldots a_{k+1} \le
\left\{\begin{array}{lll} 
q^{(k^2-1)/2k - \delta} &\text{if}\ 2 \le k \le 4, \\ 
q^{5/2 - \delta} &\text{if}\   k = 5, \\ 
q^{2(k^2-1)/3(k+2) - \delta}  &\text{if}\ k \ge 6,\end{array}\right.
\end{equation} 
for any fixed $\delta>0$,
provided that $q$ is large enough.
In fact we consider a more general case 
when $n_1, \ldots, n_k$  and $\inv{n_1 \ldots n_k}$
belong to a certain box inside of the cube $\T_{k+1}$,
where 
$$
\T_s=(\R/\Z)^s=[0, 1)^s 
$$ 
is the $k$-dimensional unit torus.

The question about the distribution of elements of $\Nab$
in various regions of $\T_k$
has  also been studied in~\cite{ASZ}.  
For an arbitrary region
$\varOmega\subseteq \T_k$  
we  denote by $\NAB$, 
the set of vectors $\vec{n} \in \Nab$
which belong to the dilated region $q\varOmega$. 
Let $\lambda(\varOmega)$ denote the Lebesgue measure of $\varOmega$.
It has been show in~\cite{ASZ} that 
for   any  fixed $k$ and region $\varOmega \subseteq \T_k$,
with piecewise smooth boundary,  
\begin{equation}
\label{eq:Omega}
\# \NAB= \lambda(\varOmega) \frac{\varphi(q)^k}{a_1 \ldots a_{k+1}}
+ O(q^{k-1/2(k+1) + o(1)})
\end{equation} 
holds uniformly over all 
vectors $\vec{a} \in \cU_q^{k+1}$ and  $\vec{b} \in
\Z^{k+1}$.

Here we show that using some bounds from~\cite{ASZ}
in a combination with some results of M.~Laczkovich~\cite{Lac}  and of
H.~Niederreiter and J.~M.~Wills~\cite{NiWi}, from the theory of uniformly 
distributed sequences leads to a better  error term 
in the asymptotic
formula~\eqref{eq:Omega}. 

Furthermore, we also consider a generalisation 
to the joint distribution of $n_1, \ldots, n_k$  and $\inv{n_1 \ldots n_k}$
in arbitrary regions. 
Namely given  an arbitrary region
$\varTheta\subseteq \T_{k+1}$ 
we  estimate the cardinality of  $\MAB$, which is
the set of vectors $(n_1, \ldots, n_k) \in \Nab$
for which $(n_1, \ldots, n_k,\inv{n_1 \ldots n_k})$
belongs to   the dilated region $q\varTheta$.

Finally, in the case of prime $q =p$, we show that a 
result of A.~Ayyad, T.~Cochrane and
Z.~Zheng~\cite[Theorem~2]{ACZ} leads to some improvements.

We conclude with a short discussion of possible ways to improve 
our results and of some open
problems.

 Throughout the paper, the implied constants in the symbols `$O$',
  and `$\ll$' may depend on  integer parameters $k$ and $r$ 
and a region $\varOmega \subseteq \T_k$.   
We recall that the notations $U = O(V)$ and $V \ll U$ are both 
equivalent to the assertion that the inequality $|U|\le cV$ holds for some
constant $c>0$.

\section{Preparations} 
\label{sec:prep}

\subsection{Character Sums} 

Let $\cX_q$ be the set of all $\# \cX_q = \varphi(q)$
multiplicative characters of $q$. We refer 
to~\cite{LN} for definitions
and basic properties of multiplicative characters
such as $\chi(u) = 0$ for any $\chi \in \cX_q$ if $\gcd (u,q) >1$
In particular, we recall that  for $u \in Z$,
\begin{equation}
\label{eq:ident}
\frac{1}{\varphi(q)}
 \sum_{\chi \in \cX_q}\chi\(u\)
=\left\{\begin{array}{rll}1 &\text{if}\ u
\equiv 1 \pmod q, \\ 0 &\text{otherwise,}\end{array}\right.
\end{equation} 
see~\cite[Theorem~5.4]{LN}.
We also use $\chi_0$ to denote the principal character.

In particular, we immediately see the following 
bound

\begin{lem}
 \label{lem:Sum Squares} For any   integers
$K$ and $L$ with $0\le K < K+ L\le q$, an integer $a\ge 1$ with $\gcd(a,q)=1$
and an arbitrary integer $b$, 
$$
\sum_{\chi \in \cX_q} 
\left| \sum_{\substack{K+1 \le n \le K+L \\ n \equiv b \pmod a}}
\chi(n)\right|^2 \le  \varphi(q)\(L/a + 1\) .
$$ 
\end{lem}

\begin{proof} We recall that if $\gcd(n,q)=1$
for the conjugated character $\overline\chi$ we
have $\overline\chi(n) =\chi\(\inv{n}\)$.
Therefore 
\begin{eqnarray*}
\sum_{\chi \in \cX_q} 
\left| \sum_{\substack{K+1 \le n \le K+L \\ n \equiv b \pmod a}}
\chi(n)\right|^2 & = &\sum_{\chi \in \cX_q} 
\left| \sum_{\substack{K+1 \le n \le K+L \\ n \equiv b \pmod a
\\\gcd(n,q)=1}}\chi(n)\right|^2\\
 & = &\sum_{\chi \in \cX_q} 
 \sum_{\substack{K+1 \le n \le K+L \\ n \equiv b \pmod a
\\\gcd(n,q)=1}} \sum_{\substack{K+1 \le m \le K+L \\ m \equiv b \pmod a
\\\gcd(m,q)=1}} \chi\(n\) \chi(\inv{m})\\
 & = & 
 \sum_{\substack{K+1 \le n \le K+L \\ n \equiv b \pmod a
\\\gcd(n,q)=1}} \sum_{\substack{K+1 \le m \le K+L \\ m \equiv b \pmod a
\\\gcd(m,q)=1}}\sum_{\chi \in \cX_q} \chi\(n \inv{m}\) = \varphi(q) T, 
\end{eqnarray*} 
where $T$ is the number of pairs $(n,m)$
with  
\begin{eqnarray*}
K+1  \le m,n \le K+L, & \qquad m &\equiv  n \equiv b \pmod a,\\
 \gcd(mn,q)=1,&\qquad  m &\equiv n \pmod q.
\end{eqnarray*} 
Clearly $n$ takes at most $L/a+1$ possible values and since 
$0\le K < K+ L\le q$, for each $n$ the 
value  of $m$ is uniquely defined.
Therefore $T \le L/a + 1$,  which concludes the proof. 
\end{proof}

The following result  is a combination of  the 
Polya-Vinogradov bound (for $r =1$) and Burgess
(for $r\ge2$) bounds, 
see~\cite[Theorems~12.5 and 12.6]{IwKow}. 

\begin{lem}
 \label{lem:PVB} For any positive integers
$U$ and $V \le q$, the bound
$$
\max_{\substack{\chi \in \cX_q\\ \chi \ne \chi_0}}
\left| \sum_{n = U+1}^{U+V}
\chi(n)\right|  \le V^{1 -1/r} q^{(r+1)/4r^2 + o(1)}. 
$$ 
holds with $r = 1,2,3$ for any $q$ and with arbitrary 
integer $r$ if $q=p$ is prime. 
\end{lem}

Now, using the identity 
\begin{equation}
\begin{split}
\label{eq:AP 2 Interv}
\sum_{\substack{K+1 \le n \le K+L \\ n \equiv b \pmod a}}
\chi(n) & = \sum_{K+1 \le am+b \le K+L}
\chi(am+b)\\
& = \chi(a)\sum_{(K+1 - b)/a \le m \le (K+L-b)/a }
\chi(m+b\inv{a})
\end{split}
\end{equation} 
we derive from Lemma~\ref{lem:PVB} the following useful estimate.

\begin{lem}
 \label{lem:Char Progr} For   any positive integers
$K$,  $L$ and  $a\ge L$ such that  $\gcd(a,q)=1$ and
an arbitrary integer $b$, the bound
$$
\max_{\substack{\chi \in \cX_q\\ \chi \ne \chi_0}}
\left| \sum_{\substack{K+1 \le n \le K+L \\ n \equiv b \pmod a}}
\chi(n)\right| \le  q^{(4r^2-3r +1)/4r^2 + o(1)} a^{-(r-1)/r} 
$$ 
holds with $r = 1,2,3$ for any $q$ and with arbitrary 
integer $r$ is $q=p$ is prime. 
\end{lem}

We also need to approximate the value of the sum of Lemma~\ref{lem:Char Progr} 
with the principal character $\chi_0$.

We denote by $\omega(q)$ the number of prime divisors 
of $q$. 

\begin{lem}
 \label{lem:Erat} For any positive integers
$K$ and $L$, an integer $a\ge 1$ such that $\gcd(a,q)=1$ and
an arbitrary integer $b$, 
$$ \sum_{\substack{K+1 \le n \le K+L \\ n \equiv b \pmod a}}
\chi_0(n)  =  \frac{\varphi(q) L}{aq}  + O(2^{\omega(q)}).
$$ 
\end{lem}

\begin{proof} Clearly
$$
\sum_{\substack{K+1 \le n \le K+L \\ n \equiv b \pmod a}}
\chi_0(n) = \sum_{\substack{K+1 \le am+b \le K+L\\\gcd(am +b, q)=1}}1.
$$

Using the M\"obius function $\mu (d)$ over the
divisors of $q$
to detect the co-primality condition 
and interchanging the order of summation, we obtain  
$$ \sum_{\substack{K+1 \le am+b \le K+L\\ \gcd(am +b, q)=1}} 1
=\sum_{d | q}\mu (d) \(\frac{L}{da} + O(1)\) = 
\frac{L}{a}  \sum_{d | q}\frac{\mu (d)}{d}
+O\( \sum_{d | q}|\mu (d)|\)
$$
from which the result follows immediately. 
\end{proof}

Finally, if $q=p$ then we also use  
the following bound which follows from a result of A.~Ayyad, T.~Cochrane and
Z.~Zheng~\cite[Theorem~2]{ACZ} and the identity~\eqref{eq:AP 2 Interv}. 

 \begin{lem}
 \label{lem:Sum 4} Let $q = p$ be prime. For any   integers
$K$ and $L$ with $0\le K < K+ L\le p$, an integer $a\ge 1$ with $\gcd(a,p)=1$
and an arbitrary integer $b$, 
$$
\sum_{\substack{\chi \in \cX_p\\ \chi \ne \chi_0}} 
\left| \sum_{\substack{K+1 \le n \le K+L \\ n \equiv b \pmod a}}
\chi(n)\right|^4 \ll  p(L/a + 1)^2 (\log p)^2 .
$$ 
\end{lem}

\subsection{Discrepancy}
\label{sec:discr}

For a finite set $\cF \subseteq \T_s$ of the unit $s$-dimensional 
set, we define its {\it discrepancy
with respect to a domain $\varXi \subseteq \T_s$\/}
as 
$$
\varDelta(\cF ,\varXi) = \left| \frac{ \#\{ \vec{f}\in\cA :\ \vec{f}\in
\varXi\} }{\#\cF} - \lambda(\varXi)\right|,
$$
where, as before, $\lambda$ is the Lebesgue measure on $\T_s$.

We now define  the {\it  discrepancy\/} of $\cF $ as 
$$
D(\cF ) = \sup_{\varPi  \subseteq \T_s}  \varDelta(\cF ,\varPi) ,
$$
where the supremum is taken over all boxes $\varPi = [\alpha_1,
\beta_1) \times \ldots \times [\alpha_s, \beta_s)$.

As usual, we define the  distance between a vector $\vec{u} \in \T_s$
and a set $\varGamma \subseteq \T_s $  by
$$
\dist(\vec{u},\varGamma) = \inf_{\vec{w} \in\varGamma}
\|\vec{u} - \vec{w}\|, 
$$
where, as before, $\|\vec{v}\|$ denotes the Euclidean norm of $\vec{v}$. Given
$\varepsilon >0$ and a domain  $\varXi \subseteq \T_s $ we define
the  sets
$$
\varXi_\varepsilon^{+} = \left\{ \vec{u} \in  \T_s \backslash
\varXi \ | \ \dist(\vec{u},\varXi) < \varepsilon \right\}
$$
and
$$
\varXi_\varepsilon^{-} = \left\{ \vec{u} \in \varXi \ | \
\dist(\vec{u},\T_s \backslash \varXi )  < \varepsilon  \right\} .
$$

Let $h(\varepsilon)$ be an arbitrary increasing function defined for
$\varepsilon >0$ and such that 
$$
\lim_{ \varepsilon \to 0}h(\varepsilon) = 0.
$$ 
As in~\cite{Lac,NiWi}, we define the
class $\cS_h$ of  domains  $\varXi \subseteq \T_s $ for which
$$
\lambda\(\varXi_\varepsilon^{+} \)\le h(\varepsilon)
\qquad \mbox{and}
\qquad
\lambda\(\varXi_\varepsilon^{-} \)\le h(\varepsilon) .
$$

A relation  between $D(\cF)$ and $\varDelta(\cF ,\varXi)$ for $\varXi \in
\cS_h$ is given by the following inequality of
M.~Laczkovich~\cite{Lac} (see also~\cite{NiWi}).

\begin{lem}
 \label{lem:LNW bound} For any domain  $\varXi \in \cS_h$, we have
$$
\varDelta(\cF ,\varXi) \ll h\(s^{1/2} D(\cF )^{1/s}\)  .
$$ 
\end{lem}

Finally,  the following bound, which is a partial case 
of a more general result of H.~Weyl~\cite{Weyl}
shows that if  $\varXi$ has a piecewise smooth boundary
that $\varXi \in \cS_h$ for some linear function 
$h(\varepsilon) = C\varepsilon$. 

\begin{lem}
 \label{lem:Weyl bound} For any domain  $\varXi \in \cS_h$ with
piecewise smooth boundary, we
have
$$
\lambda\(\varXi_\varepsilon^{\pm }\) = O(\varepsilon).
$$ 
\end{lem}

\section{Main Results}

\subsection{Distribution in Boxes}

Here we study $\MAB$ in the case where
$\varTheta = \varSigma$ is  a box $\varSigma \subseteq \T_{k+1}$
and in particular 
we generalise and improve the bound~\eqref{eq:Old Nab}.

We recall that we use $\|\vec{a}\|$ to denote the Euclidean norm of  
$\vec{a} \in \cU_q^{k+1}$. 
  
\begin{thm}
\label{thm:main-1}  
For $r=1,2,3$,  any  fixed $k\ge 2$, 
and a box 
$$\varSigma =   [\alpha_1,
\beta_1) \times \ldots \times [\alpha_{k+1}, \beta_{k+1})
\subseteq \T_{k+1}
$$
the bound  
\begin{eqnarray*}
\lefteqn{\# \MABB = \lambda(\varSigma) \frac{\varphi(q)^k}{a_1 \ldots
a_{k+1}}}\\
& & \qquad \qquad  + O\(\frac{\|\vec{a}\| q^{k-1 + o(1)}}{a_1 \ldots a_{k+1}}
+\frac{q^{k - (3r - 1)(k-1)/4r^2 + o(1)}}
{(a_1\ldots a_{k+1})^{1- (k + r - 1)/r(k+1)}} \)
\end{eqnarray*} 
holds uniformly over all 
vectors $\vec{a} \in \cU_q^{k+1}$ and  $\vec{b} \in
\Z^{k+1}$.
\end{thm}

\begin{proof} We see from~\eqref{eq:ident}  that  
\begin{eqnarray*}
\lefteqn{\# \MABB = 
\sum_{\substack{{\alpha_1q} \le n_1 <  {\beta_1 q}\\
n_1 \equiv b_1 \pmod {a_1} }}\ldots
\sum_{\substack{ {\alpha_{k+1}q} \le n_{k+1} <  {\beta_{k+1} q}\\
n_{k+1} \equiv b_{k+1}\pmod {a_{k+1}}}}} \\
& &\qquad \qquad \qquad \qquad \qquad \qquad \qquad \qquad 
\frac{1}{\varphi(q)}
 \sum_{\chi \in \cX_q}\chi\(n_1\ldots n_{k+1} \).
\end{eqnarray*} 

We now change the order of summation and  note that by 
Lemma~\ref{lem:Erat} the term corresponding to 
the principal character $\chi =\chi_0$ is equal to 
\begin{eqnarray*}
\lefteqn{\frac{1}{\varphi(q)} \sum_{\substack{ {\alpha_1q} \le n_1 <  {\beta_1
q}\\ n_1 \equiv b_1 \pmod {a_1}  \\ \gcd(n_1,q) =1}}\ldots
\sum_{\substack{ {\alpha_{k+1}q} \le n_{k+1} <  {\beta_{k+1} q}\\
n_{k+1} \equiv b_{k+1}\pmod {a_{k+1}} \\ \gcd(n_{k+1},q) =1}} 1}\\
&  & \qquad \qquad 
=\frac{1}{\varphi(q)} \prod_{\nu=1}^{k+1} \(\frac{(\beta_\nu
-\alpha_\nu) \varphi(q)}{a_\nu}  + O\(2^{\omega(q)}\)\)\\
&  & \qquad \qquad 
 = \lambda(\varSigma) \frac{\varphi(q)^k}{a_1 \ldots a_{k+1}}
+ O\(\frac{\|\vec{a}\| \varphi(q)^{k-1}}{a_1 \ldots a_{k+1}} 2^{k
\omega(q)}\)\\
&  & \qquad \qquad 
= \lambda(\varSigma) \frac{\varphi(q)^k}{a_1 \ldots a_{k+1}}
+ O\(\frac{\|\vec{a}\| q^{k-1 + o(1)}}{a_1 \ldots a_{k+1}} \),
\end{eqnarray*} 
since  
$$
\omega(q) \ll \frac{\log q}{\log \log q},
$$
see~\cite[Section~I.5.3]{Ten}, and the bound~\eqref{eq:phi}.
Hence,
\begin{equation}
\label{eq:Bound M and R}
\# \MABB = \lambda(\varSigma) \frac{\varphi(q)^k}{a_1 \ldots a_{k+1}}
+ O\(\frac{\|\vec{a}\| q^{k-1 + o(1)}}{a_1 \ldots a_{k+1}} + R\),
\end{equation}
where 
\begin{eqnarray*}
R & = &\frac{1}{\varphi(q)}  \sum_{\substack{\chi \in
\cX_q\\\chi\ne \chi_0}}  
\sum_{\substack{ {\alpha_1q} \le n_1 <  {\beta_1 q}\\
n_1 \equiv b_1 \pmod {a_1} }}\ldots
\sum_{\substack{ {\alpha_{k+1}q} \le n_{k+1} <  {\beta_{k+1} q}\\
n_{k+1} \equiv b_{k+1}\pmod {a_{k+1}}}}\chi\(n_1\ldots n_{k+1}
\)  \\
 & = & \frac{1}{\varphi(q)}  \sum_{\substack{\chi \in
\cX_q\\\chi\ne \chi_0}}
\prod_{\nu=1}^{k+1}  
\sum_{\substack{ {\alpha_\nu q} \le n_\nu <  {\beta_\nu q}\\
n_\nu \equiv b_\nu \pmod {a_\nu} }} \chi\(n_\nu\).
\end{eqnarray*}

Thus using the H{\"o}lder inequality we obtain 
\begin{equation}
\label{eq:Bound R q}
 R^{k+1} \le \frac{1}{\varphi(q)^{k+1}}  \prod_{\nu=1}^{k+1}
\sum_{\substack{\chi \in \cX_q\\\chi\ne \chi_0}}
\left|\sum_{\substack{ {\alpha_\nu q} \le n_\nu <  {\beta_\nu q}\\
n_\nu \equiv b_\nu \pmod {a_\nu} }} \chi\(n_\nu\)\right|^{k+1}.
\end{equation}

Since the bound is trivial for $\|\vec{a}\|\ge q$,
we can assume that 
$$
\max\{a_1,  \ldots,  a_{k+1}\} < q. 
$$

Applying Lemma~\ref{lem:Char Progr} to the $(k-1)$th power of 
the character sums  for each
$\nu=1, \ldots, k+1$, and then extending the summation over all
characters $\chi \in \cX$,  we obtain 
that for $r=1,2,3$
\begin{eqnarray*}
\lefteqn{
\sum_{\substack{\chi \in \cX_q\\\chi\ne \chi_0}}
\left|\sum_{\substack{ {\alpha_\nu q} \le n_\nu <  {\beta_\nu q}\\
n_\nu \equiv b_\nu \pmod {a_\nu} }} \chi\(n_\nu\)\right|^{k+1}}\\
& & \qquad \le \(q^{(4r^2-3r +1)/4r^2 + o(1)} a_\nu^{-(r-1)/r}\)^{k-1}  
\sum_{\chi \in \cX_q}
\left|\sum_{\substack{ {\alpha_\nu q} \le n_\nu <  {\beta_\nu q}\\
n_\nu \equiv b_\nu \pmod {a_\nu} }} \chi\(n_\nu\)\right|^{2}.
\end{eqnarray*}
We now use Lemma~\ref{lem:Sum Squares}, which implies
\begin{eqnarray*}
\sum_{\chi \in \cX_q}
\left|\sum_{\substack{ {\alpha_\nu q} \le n_\nu <  {\beta_\nu q}\\
n_\nu \equiv b_\nu \pmod {a_\nu} }} \chi\(n_\nu\)\right|^{k+1}& \le &
\(q^{(4r^2-3r +1)/4r^2 + o(1)} a_\nu^{-(r-1)/r}\)^{k-1}  q^2  a_\nu^{-1}\\
 & \le &
q^{(k+1) - (3r - 1)(k-1)/4r^2 + o(1)} a_\nu^{-(rk -k + 1)/r} .
\end{eqnarray*}
Substituting this bound in~\eqref{eq:Bound R q} and using~\eqref{eq:phi},  we
obtain 
\begin{equation*}
\begin{split}
R & \le 
q^{k - (3r - 1)(k-1)/4r^2 + o(1)}
\(\prod_{\nu=1}^{k+1}a_\nu\)^{- (rk - k + 1)/r(k+1)},
\end{split}
\end{equation*}
which together with~\eqref{eq:Bound M and R} completes the proof.
\end{proof}

In particular, taking $r=3$  we see that the bound of Theorem~\ref{thm:main-1} 
implies that for any fixed
$k$ and 
$\delta >0$ there exists $\eta>0$ such that  under
the conditions~\eqref{eq:New Threshold}
we have 
$$
\# \MABB= \(\lambda(\varSigma) + O(q^{-\eta})\) \frac{\varphi(q)^k}{a_1
\ldots a_{k+1}}  .
$$
Moreover, using the trivial bounds
$$
\|\vec{a}\| \le a_1 \ldots a_{k+1} \qquad \text{and}
\qquad (a_1\ldots
a_{k+1})^{1/(k+1)} \ge 1, 
$$
we derive from Theorem~\ref{thm:main-1}  that 
$$
\# \MABB= \lambda(\varSigma) \frac{\varphi(q)^k}{a_1 \ldots a_{k+1}}
+ O\(q^{k-1 + o(1)} + q^{(k+1)/2 +o(1)} \).
$$
Finally, taking   $\varSigma =   \T_{k+1}$  and $r=1$ in
Theorem~\ref{thm:main-1},  
 we obtain~\eqref{eq:New Nab}.

\subsection{Distribution in General Regions}

Here we give an improvement and 
generalisation of  the asymptotic
formula~\eqref{eq:Omega}.

\begin{thm}
\label{thm:main-2}  
For  $r=1,2,3$,  any  fixed $k\ge 2$ 
and region $\varTheta \subseteq \T_{k+1}$ 
with piecewise smooth boundary,   
\begin{eqnarray*}
\lefteqn{\# \MAB = \lambda(\varTheta) \frac{\varphi(q)^k}{a_1 \ldots a_{k+1}}}\\
& & \qquad \qquad +~O\(\frac{\|\vec{a}\|^{1/(k+1)}q^{k-1/(k+1) +o(1)}}{a_1
\ldots a_{k+1}} +\frac{q^{k-(3r - 1)(k-1)/4r^2(k+1)  +
o(1)}}{(a_1\ldots a_{k+1})^{1- (k + r - 1)/r(k+1)^2}} \)
\end{eqnarray*} 
holds uniformly over all 
vectors $\vec{a} \in \cU_q^{k+1}$ and  $\vec{b} \in
\Z^{k+1}$.
\end{thm}

\begin{proof} It follows from Theorem~\ref{thm:main-1}
and the estimates~\eqref{eq:phi}  and~\eqref{eq:New Nab}
that we have  the bound 
$$
D(\cA(\vec{a},\vec{b},q)) \ll \|\vec{a}\|q^{-1 +o(1)}  +
(a_1\ldots a_{k+1})^{(r+k -1)/r(k+1)}q^{-(3r - 1)(k-1)/4r^2 + o(1)}
$$
on the box discrepancy of the set
$$
 \cA(\vec{a},\vec{b},q) = 
\left\{\(\frac{n_1}{q}, \ldots, \frac{n_k}{q}, 
\frac{\inv{n_1 \ldots n_k}}{q}\) \ :
(n_1, \ldots, n_k) \in \Nab\right\}.
$$

Therefore, by Lemma~\ref{lem:LNW bound} and Lemma~\ref{lem:Weyl bound}
we conclude  that the  discrepancy of $\cA(\vec{a},\vec{b},q)$
with respect to $\varTheta$ satisfies
\begin{eqnarray*}
\lefteqn{
\varDelta(\cA(\vec{a},\vec{b},q),\varTheta) \ll \|\vec{a}\|^{1/(k+1)}q^{-1/(k+1)
+o(1)}  }\\
& & \qquad\qquad\qquad \qquad
+~(a_1\ldots a_{k+1})^{k/r(k+1)^2}q^{-(3r - 1)(k-1)/4r^2(k+1) +o(1)}
\end{eqnarray*}
which is equivalent to the desired result. 
\end{proof} 

Certainly applying Theorem~\ref{thm:main-2} with 
$\varTheta = \varOmega \times [0,1)$ 
where  $\varOmega \subseteq \T_k$ one immediately
obtains an asymptotic formula for $\NAB$ (which is already
stronger than~\eqref{eq:Omega}. However since the problem 
of estimating $\NAB$ is of lower dimension
($k$ instead of $k+1$) one obtains a slightly 
stronger bound in this case.

\begin{thm}
\label{thm:main-3}  
For  $r=1,2,3$,    any  fixed $k\ge $ 
and region $\varOmega\subseteq \T_k$ 
with piecewise smooth boundary,   
\begin{eqnarray*}
\# \NAB & = & \lambda(\varOmega) \frac{\varphi(q)^k}{a_1 \ldots a_{k+1}}\\
& & \quad + O\(\frac{\|\vec{a}\|^{1/k}q^{k-1/k +o(1)}}{a_1 \ldots
a_{k+1}} +\frac{q^{k-(3r - 1)(k-1)/4r^2k +
o(1)}}{(a_1\ldots a_{k+1})^{1- (k + r - 1)/rk(k+1)}} \)
\end{eqnarray*} 
holds uniformly over all 
vectors $\vec{a} \in \cU_q^{k+1}$ and  $\vec{b} \in
\Z^{k+1}$.
\end{thm}

\begin{proof} Taking the box $\varSigma = \varPi \times [0,1)$,
 where a box $\varPi \subseteq
\T_k$ we see that it follows from Theorem~\ref{thm:main-1}
and the estimates~\eqref{eq:phi}  and~\eqref{eq:New Nab}
that we have the bound 
$$
D(\cB(\vec{a},\vec{b},q)) \ll \|\vec{a}\|q^{-1 +o(1)}  +
(a_1\ldots a_{k+1})^{ (r+k -1)/r(k+1)}q^{-(3r - 1)(k-1)/4r^2 + o(1)}
$$
on the box discrepancy of the set
$$
 \cB(\vec{a},\vec{b},q) = 
\left\{\(\frac{n_1}{q}, \ldots, \frac{n_k}{q} \) \ :
(n_1, \ldots, n_k) \in \Nab\right\}.
$$
that is,  of the same strength as that 
for the set $\cA(\vec{a},\vec{b},q)$ in the proof 
of  Theorem~\ref{thm:main-2}.

Therefore, by Lemma~\ref{lem:LNW bound} and Lemma~\ref{lem:Weyl bound}
we conclude  that the  discrepancy of $\cB(\vec{a},\vec{b},q)$
with respect to $\varOmega$ satisfies
\begin{eqnarray*}
\lefteqn{
\varDelta(\cB(\vec{a},\vec{b},q),\varOmega) }\\
& & \qquad \ll \|\vec{a}\|^{1/k}q^{-1/k
+o(1)}  + (a_1\ldots a_{k+1})^{(r+k -1)/rk(k+1)}q^{-(3r - 1)(k-1)/4r^2k + o(1)}
\end{eqnarray*}
which is equivalent to the desired result. 
\end{proof}

We remark that although in the case of $k=2$ and fixed $\vec{a}$, 
Theorem~\ref{thm:main-1} (with the optimal choice of $r=1$)
is equivalent to~\eqref{eq:Old Nab}, 
the bound of Theorem~\ref{thm:main-3}
still improves~\eqref{eq:Omega} due to our use of Lemma~\ref{lem:LNW bound}
instead of the arguments from~\cite{ASZ}.

\subsection{Some Improvements for Prime $q=p$}

Here we show that if $q=p$ is prime and $k\ge 3$
then using Lemma~\ref{lem:Sum 4} instead of 
Lemma~\ref{lem:Sum Squares} leads to a stronger bounds 
wit respect to the product $a_1\ldots a_{k+1}$.

\begin{thm}
\label{thm:main-4}  Let $q=p$ be prime.
For   any  fixed integer $r \ge 1$ and $k\ge 3$, 
and a box 
$$\varSigma =   [\alpha_1,
\beta_1) \times \ldots \times [\alpha_{k+1}, \beta_{k+1})
\subseteq \T_{k+1}
$$
the bound 
\begin{eqnarray*}
\lefteqn{
\# \MABBp= \lambda(\varSigma) \frac{p^k}{a_1 \ldots a_{k+1}}}\\
& & \quad +~O\(\frac{\|\vec{a}\| p^{k-1 + o(1)}}{a_1 \ldots a_{k+1}} 
+\frac{p^{k - (3r - 1)(k-3)/4r^2 + o(1)}}
{(a_1\ldots a_{k+1})^{1-(k + 2r -3)/r(k+1)}} \)
\end{eqnarray*}
holds uniformly over all 
vectors $\vec{a} \in \cU_q^{k+1}$ and  $\vec{b} \in
\Z^{k+1}$.
\end{thm}

\begin{proof} We proceed as in the proof of Theorem~\ref{thm:main-1}
and also note that one can change $(p-1)^k$ to $p^k$ 
in~\eqref{eq:Bound M and R}
without changing the error term. In particular, we still 
have~\eqref{eq:Bound R q}. Again as in the proof of Theorem~\ref{thm:main-1}
we apply Lemma~\ref{lem:Char Progr}, 
however this time to the $(k-3)$th power of 
the character sums  for each
$\nu=1, \ldots, k+1$ (and this time we do not extend   the summation over all
characters $\chi \in \cX$),  we obtain that for any integer $r\ge 1$
\begin{eqnarray*}
\lefteqn{
\sum_{\substack{\chi \in \cX_p\\\chi\ne \chi_0}}
\left|\sum_{\substack{ {\alpha_\nu p} \le n_\nu <  {\beta_\nu p}\\
n_\nu \equiv b_\nu \pmod {a_\nu} }} \chi\(n_\nu\)\right|^{k+1}}\\
& & \qquad \le \(p^{(4r^2-3r +1)/4r^2 + o(1)} a_\nu^{-(r-1)/r}\)^{k-3}  
\sum_{\substack{\chi \in \cX_p\\\chi\ne \chi_0}}
\left|\sum_{\substack{ {\alpha_\nu p} \le n_\nu <  {\beta_\nu p}\\
n_\nu \equiv b_\nu \pmod {a_\nu} }} \chi\(n_\nu\)\right|^{4}.
\end{eqnarray*}
We now use Lemma~\ref{lem:Sum 4}, which implies
\begin{eqnarray*}
\sum_{\substack{\chi \in \cX_p\\\chi\ne \chi_0}}
\left|\sum_{\substack{ {\alpha_\nu p} \le n_\nu <  {\beta_\nu p}\\
n_\nu \equiv b_\nu \pmod {a_\nu} }} \chi\(n_\nu\)\right|^{k+1}& \le &
\(p^{(4r^2-3r +1)/4r^2 + o(1)} a_\nu^{-(r-1)/r}\)^{k-3}  p^{4} 
a_\nu^{-2}\\
 & \le &
p^{k+1 - (3r - 1)(k-3)/4r^2 + o(1)} a_\nu^{-(rk - k -r +3)/r(k+1)}. 
\end{eqnarray*}
Substituting this bound in~\eqref{eq:Bound R q} and using~\eqref{eq:phi},  we
obtain 
\begin{equation*}
\begin{split}
R & \le 
p^{k - (3r - 3)(k-1)/4r^2 + o(1)}
\(\prod_{\nu=1}^{k+1}a_\nu\)^{-(rk - k -r +3)/r(k+1)},
\end{split}
\end{equation*}
which together with~\eqref{eq:Bound M and R} completes the proof.
\end{proof}

In particular we see that the bound of Theorem~\ref{thm:main-4} 
taken with $r=1$ 
implies that for any fixed
$k$ and 
$\delta >0$ there exists $\eta>0$ such that  under
the conditions
$$
\|\vec{a}\| \le p^{1 - \delta} \qquad 
\text{and} \qquad   a_1 \ldots a_{k+1} \le p^{(k+1)/2 - \delta}
$$
we have 
$$
\# \MABBp= \(\lambda(\varSigma) + O(p^{-\eta})\) \frac{p^k}{a_1
\ldots a_{k+1}}  .
$$
However,   taking a sufficiently large $r$ we see from 
Theorem~\ref{thm:main-4}  that for any $\delta >0$ there exists
$K_0$ and $\eta>0$ such that for $k \ge K_0$ the above bound
holds under the condition 
$$
\|\vec{a}\| \le p^{1 - \delta} \qquad 
\text{and} \qquad   a_1 \ldots a_{k+1} \le p^{(3/4 - \delta)k}
$$

We also have analogues of Theorems~\ref{thm:main-2} and~\ref{thm:main-3}. 

\begin{thm}
\label{thm:main-5}  Let $q=p$ be prime.
For   any  fixed integer $r \ge 1$ and $k\ge 3$, 
and region $\varTheta \subseteq \T_{k+1}$ 
with piecewise smooth boundary,   
\begin{eqnarray*}
\lefteqn{\# \MABp = \lambda(\varTheta) \frac{p^k}{a_1 \ldots a_{k+1}}}\\
& & \qquad \quad +~O\(\frac{\|\vec{a}\|^{1/(k+1)}p^{k-1/(k+1) +o(1)}}{a_1
\ldots a_{k+1}} +\frac{p^{k - (3r - 1)(k-3)/4r^2(k+1) + o(1)}}
{(a_1\ldots a_{k+1})^{1 - (k+2r-3)/r(k+1)^2}} \)
\end{eqnarray*} 
holds uniformly over all 
vectors $\vec{a} \in \cU_q^{k+1}$ and  $\vec{b} \in
\Z^{k+1}$.
\end{thm}

\begin{thm}
\label{thm:main-6}  Let $q=p$ be prime.
For   any  fixed integer $r \ge 1$ and $k\ge 3$, 
and region $\varOmega\subseteq \T_k$ 
with piecewise smooth boundary,   
\begin{eqnarray*}
\lefteqn{
\# \NABp   =   \lambda(\varOmega) \frac{p^k}{a_1 \ldots a_{k+1}}}\\
& & \qquad  +~O\(\frac{\|\vec{a}\|^{1/k}p^{k-1/k +o(1)}}{a_1 \ldots
a_{k+1}} +\frac{p^{k - (3r - 1)(k-3)/4r^2k + o(1)}}
{(a_1\ldots a_{k+1})^{1 - (k+2r-3)/rk(k+1)}}  \)
\end{eqnarray*} 
holds uniformly over all 
vectors $\vec{a} \in \cU_q^{k+1}$ and  $\vec{b} \in
\Z^{k+1}$.
\end{thm}

\section{Concluding Remarks}

\subsection{Further Improvements}

Clearly, if some of $a_1, \ldots, a_{k+1}$ are 
of different order magnitude, then in 
the proofs of Theorems~\ref{thm:main-1} and~\ref{thm:main-4} one can use
Lemma~\ref{lem:PVB} with various values of $r$ for each $\nu$ which may
lead to stronger bounds. However it seems that the optimal strategy of
applying these results heavily depends on various relations
between the sizes of $a_1,\ldots, a_{k+1}$ and $q$.

We believe that there are several further  possibilities 
of improving Theorem~\ref{thm:main-1} and in particular 
improving the threshold~\eqref{eq:New Threshold}.
Certainly there should be a variant of the result of
 A.~Ayyad, T.~Cochrane and
Z.~Zheng~\cite[Theorem~2]{ACZ}, given in Lemma~\ref{lem:Sum 4},
which holds 
for arbitrary composite moduli $q$ (see also~\cite{CZ}).  
In fact, J.~B.~Friedlander and H.~Iwaniec~\cite{FrIw}
give such a bound, but only for special intervals 
(starting at the origin).  Certainly, obtaining such a 
general result is of independent interest.
Furthermore, the technique used   by
M.~Z.~Garaev~\cite{Gar} can probably be useful as well. 

We also note that by a result of 
H.~Niederreiter and J.~M.~Wills~\cite{NiWi} for the class of convex sets
$\varTheta$ and $\varOmega$ the implied constants do   not depend 
on the set (see also~\cite[Theorem~1.12]{DrTi} 
and~\cite[Theorem~1.6, Chapter~2]{KuNi}). 

\subsection{Some Open Problems}

It is certainly interesting to study various geometric
properties of the set $\Nab$. For example,
let 
$$\Hab = \max_{(n_1, \ldots n_k) \in \Nab} 
\min_{1 \le i \le k}  \left|n_i - \inv{n_1 \ldots n_k}\right|.
$$
For $k = 1$, $\vec{a} = (1,1)$, $\vec{b} = (0,0)$, the value
of
$$
H(q) = \max_{n \in \cU_q} 
 \left|n  - \inv{n}\right|
$$
has been studied in~\cite{FKSY,KhShp}. In particular,
it has been shown in~\cite{KhShp} that $H(q) = q + O\(q^{3/4
+ o(1)}\)$. It has also been shown  in~\cite{FKSY} that
$H(q)$ is influenced by the distribution of divisors 
of $qs -1$ for small values of $s$, and thus some 
lower bounds on $H(q)$ have been derived. It would be 
interesting to find out whether the behavior of $\Hab$
is also influenced by some arithmetic properties
of the modulus $q$.

Finally, one can also study various geometric 
properties of the convex closure of
$\Nab$. For example, 
for $k = 1$, $\vec{a} = (1,1)$, $\vec{b} = (0,0)$,
that is, for the set 
$$
\cN(q) = \left\{\(n,\inv{n}\) \ : \ n \in \cU_q\right\}
$$
some lower and upper bounds on 
the number of vertices $V(q)$ of its
convex closure  have been given in~\cite{KSY}. 
These bounds as well as some numerical calculations
suggest that the convex closure of $\cN(q)$ does not behave
as the convex closure of a random set, but rather is
affected by the arithmetic structure of $q-1$ (and probably
of $qs-1$ for small integers $s$).  It would be interesting to
see whether the same effect appears in the behaviour
of the convex closure of $\Nab$ for larger  values of $k$
and ``generic'' vectors $\vec{a}$ and $\vec{b}$.



\begin{thebibliography}{100}

\bibitem{ASZ} E.  Alkan, F.~Stan and A. Zaharescu, 
`Lehmer $k$-tuples', 
{\it Proc. Amer. Math. Soc.\/}, {\bf 134} (2006), 2807--2815. 



\bibitem{ACZ} A. Ayyad, T. Cochrane and Z. Zheng, 
`The congruence $x_1x_2 \equiv x_3x_4 \pmod p$, the equation
$x_1x_2 = x_3x_4$ and the mean value of character sums', 
{\it J. Number Theory\/}, {\bf 59} (1996), 398--413.

%
%

\bibitem{CoZa} C. Cobeli and A. Zaharescu, `Generalization of a 
problem of Lehmer', 
{\it Manuscr. Math.\/}, {\bf 104} (2001),  301--307.

\bibitem{CZ}  T. Cochrane and Z. Zheng, 
 `High order moments of character sums', 
{\it Proc. Amer. Math. Soc.\/}, {\bf 126} (1998),   951--956.

\bibitem{DrTi} M. Drmota and R. Tichy, 
{\it Sequences, discrepancies and applications\/},
Springer-Verlag, Berlin, 1997.

\bibitem {FKSY} K. Ford, M. R. Khan, 
I. E. Shparlinski and C. L. Yankov,
`On the maximal difference between an element and its inverse in
residue rings', {\it Proc. Amer. Math. Soc.\/}, {\bf 133} (2005),
3463--3468.

\bibitem{FrIw}
J.~B.~Friedlander and H.~Iwaniec,
`The divisor problem for arithmetic progressions',
{\it Acta Arith.\/}, {\bf 45} (1985), 273--277.


\bibitem{Gar} M.~Z.~Garaev, `Character sums in short 
intervals and the multiplication table modulo a large prime',
{\it Monatsh. Math.\/}, {\bf  148}  (2006), 127--138.

\bibitem{Guy} R. K. Guy, {\it Unsolved problems in number theory\/},
Springer-Verlag, New York, 1994.

\bibitem{IwKow} H. Iwaniec and E. Kowalski,  {\it Analytic number
theory\/},  Amer. Math. Soc., Providence, RI, 2004.

\bibitem{KhShp} M. R. Khan and I. E. Shparlinski, `On the maximal
difference between an element and its inverse modulo $n$',
{\it Period. Math. Hung.\/}, {\bf 47} (2003), 111--117.

\bibitem {KSY}   M. R. Khan, I. E. Shparlinski and 
C. L. Yankov,
`On the convex closure of   the graph of  modular 
inversions', {\it Preprint\/}, 2006.

\bibitem{KuNi}
L. Kuipers and H. Niederreiter, {\it Uniform
distribution of sequences\/}, Wiley-Interscience Publ., 1974.

\bibitem{Lac} M. Laczkovich,
`Discrepancy estimates for sets with small
boundary', {\it Studia Sci. Math. Hungar.\/}, {\bf 30} (1995),
105--109.

\bibitem{LN} R. Lidl and H. Niederreiter, {\it Finite fields\/},
Cambridge University Press, Cambridge, 1997. 

\bibitem{NiWi} H. Niederreiter and J. M. Wills,
`Diskrepanz und Distanz von Massen bezuglich
konvexer und Jordanscher Mengen',
{\it Math. Z.\/}, {\bf 144} (1975),   125--134. 

\bibitem{Ten} G. Tenenbaum, {\it Introduction to
analytic and probabilistic
number theory\/}, Cambridge Univ. Press, 1995.

\bibitem{Weyl} H. Weyl, `On the volume of tubes',  
{\it Amer. J. Math.\/}, {\bf 61} (1939), 461--472.

\bibitem {Zha1}  W. Zhang, 
`On a problem of D.~H.~Lehmer and its generalization',  
{\it Compos. Math.\/}, {\bf 86} (1993),  307--316. 

\bibitem {Zha2}  W. Zhang, 
`On a problem of D.~H.~Lehmer and its generalization, II',  
{\it Compos. Math.\/}, {\bf 91} (1994),   47--56. 

\bibitem {Zha3}  W. Zhang, 
`On the difference between a D. H. Lehmer number and its 
inverse modulo  $q$',  
{\it Acta Arith.\/}, {\bf 68} (1994),   255--263.
\end{thebibliography}
\end{document}